%% file: appendix.tex
\documentclass[11pt,a4paper]{article}
\usepackage{amsfonts,amsmath,amsthm,epsf,graphics,verbatim,amssymb,amscd,eucal,bbm}

\newtheorem{theorem}{Theorem}[section]

\newtheorem{lemma}[theorem]{Lemma}
\theoremstyle{definition}

\theoremstyle{remark}

\def\Blackboardfont{\mathbb}

\newcommand{\moins}{ {\setminus} }

\def\Z{{\Blackboardfont Z}}

\def\N{{\Blackboardfont N}}
\def\R{{\Blackboardfont R}}
\def\T{{\Blackboardfont T}}

\def\cB{{\mathcal B}}

\def\cS{{\mathcal S}}

\def\iff{\Longleftrightarrow}
\def\eref#1{(\ref{#1})}

\begin{document}
\sloppy

\title{\bf Appendix to the paper ``Random walks on free products of cyclic groups''}

\author{Jean {\sc Mairesse}
\thanks{LIAFA, CNRS-Universit\'e Paris 7, case
    7014, 2, place Jussieu, 75251 Paris Cedex 05, France. E-mail: {\tt Jean.Mairesse@liafa.jussieu.fr}} 
    \ and Fr\'ed\'eric {\sc Math\'eus} 
\thanks{LMAM, Universit\'e de Bretagne-Sud, Campus de Tohannic,  BP
  573, 56017 Vannes, France. E-mail: {\tt Frederic.Matheus@univ-ubs.fr}}}

\maketitle

%% \centerline{\large{Draft version: Please do not circulate}}

\begin{abstract}
This paper is an appendix to the paper ``Random walks on free products
of cyclic groups'' by J. Mairesse and F. Math\'eus. 
It contains the details of the computations and the proofs of the
results concerning the examples treated there. 
\end{abstract}

\smallskip


\smallskip

\textsl{AMS classification (2000):} Primary 60J10, 60B15, 60J22, 65C40;
Secondary 28D20; 37M25.

%\newpage
%\tableofcontents
%\baselineskip 12pt

\section{Introduction}
\label{se-intro} 

Let $G$ be a free product of a finite family of finite groups, and let
the set of generators $\Sigma$ be the union of the finite groups. 
We consider a transient nearest neighbor random walk on $(G,\Sigma)$.  
The harmonic measure is Markovian with a special combinatorial
structure that can be completely described via a finite set of
polynomial equations (see \cite{MaMa} for a proof and references). 
In several simple cases of interest, the polynomial equations
can be explicitly solved, to get closed form formulas for the drift,
the entropy, the probability of ever hitting an element, or the minimal
positive harmonic functions of the walk. To illustrate, we give in
\cite{MaMa} explicit formulas for the drift in the following cases:
the modular group $\Z/2\Z\star \Z/3\Z$, $\Z/3\Z\star \Z/3\Z$,
$\Z/k\Z\star \Z/k\Z$ and the Hecke groups  
$\Z/2\Z\star \Z/k\Z$. The present paper is an appendix to \cite{MaMa}.
The purpose is to give the details of the computations and the proofs of the
results concerning these examples. Vershik's notion of extremal generators
is also investigated in \cite{MaMa}, and several examples and counter-examples
are given. Again, the details of the computations are to be found in
this appendix. 

This appendix is self-contained in the following sense: the examples
treated and the results proved are specified. On the other
hand, we make a 
constant use of notations, notions or results from \cite{MaMa} which are not
redefined or reproved. Hence, it is probably not a good idea to read this
appendix independently of the parent paper. 

\section{Explicit Drift Computations}
\label{se-expli}

We give here the details of the computations and the complete
proofs of the results stated in \cite[\S 4.2-4.5]{MaMa}.

\subsection{Random walks on $\Z/2\Z\star \Z/3\Z$}\label{sse-z2z3}

The group
$\Z/2\Z\star \Z/3\Z$  is isomorphic to the modular group PSL$(2,\Z)$,
i.e. the group of 2x2 matrices with integer entries and determinant
$1$, quotiented by $\pm$Id. Let $a$ and $b$ be
the respective generators of $\Z/2\Z$ and $\Z/3\Z$. 
A possible faithful representation of the
group is 
\begin{equation}\label{eq-modular}
a=\left[ \begin{array}{cc} 0 & -1 \\ 1 & 0 \end{array}\right], \qquad 
b=\left[ \begin{array}{cc} 1 & -1 \\ 1 & 0 \end{array}\right] \:.
\end{equation}

%%Quoting \cite[Chapter II.B]{harp}: ``The modular group is one of the
%%most important groups in mathematics''. 

%% Since this is the first case treated, we detail the
%% computations much more that in the forthcoming examples. 

\medskip

\begin{figure}[ht]
\[ \epsfxsize=220pt \epsfbox{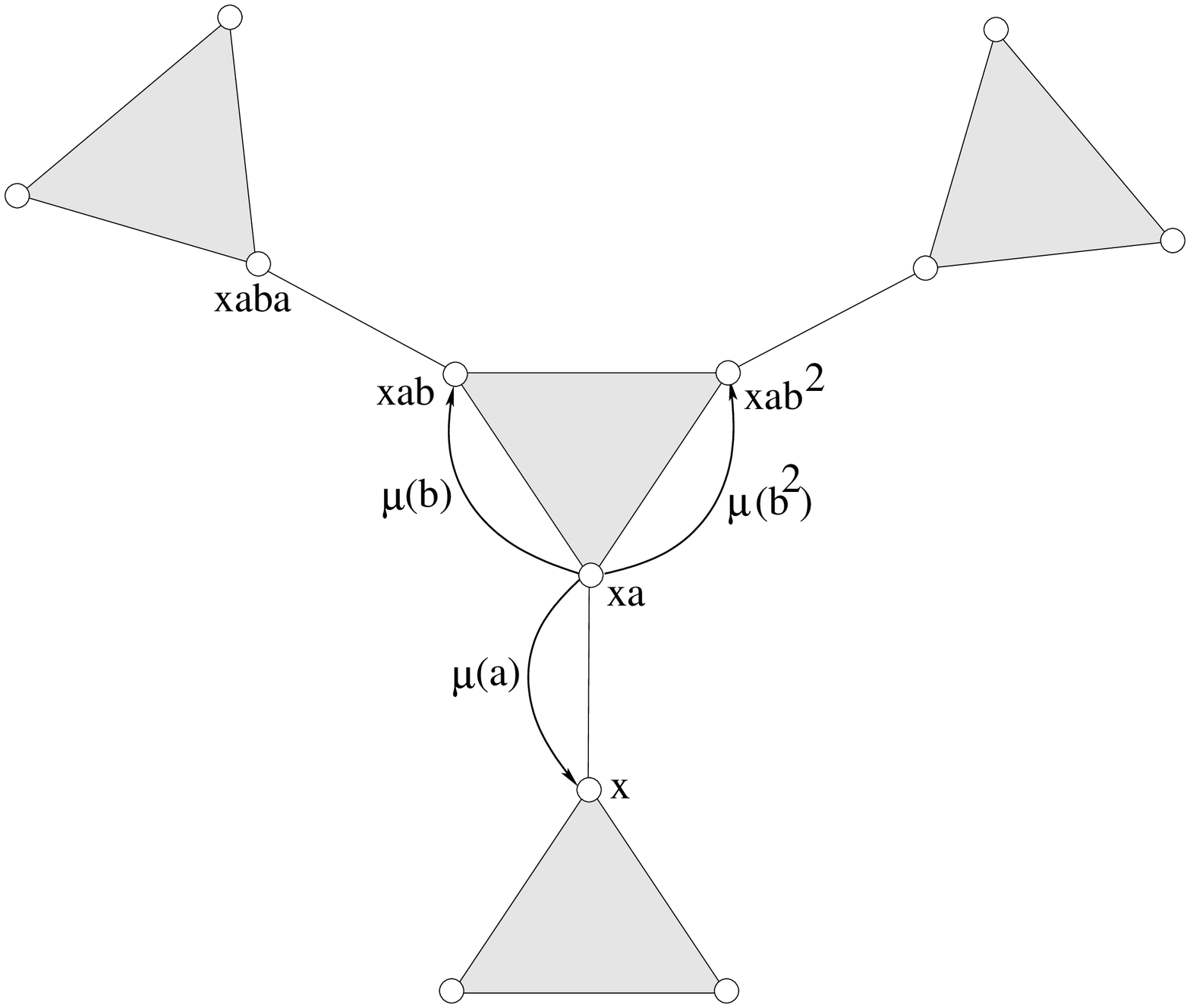} \]
\caption{A nearest neighbor random walk on $\Z/2\Z\star \Z/3\Z$.}
\label{fi-z2z3} 
\end{figure}

Consider a general nearest neighbor random walk $(\Z/2\Z\star
\Z/3\Z,\mu)$. Set $\mu(a)=1-p-q, \mu(b)=p,\mu(b^2)=q$. In
Figure \ref{fi-z2z3}, we have represented the Cayley graph of
$\Z/2\Z\star \Z/3\Z$ and the one-step transitions of the random walk. 

Consider first the case of a random walk satisfying:
$\mu(a)=1-2p,\mu(b)=\mu(b^2)=p, p\in (0,1/2)$.
Here, we are in the situation described at the beginning of \cite[\S 4.1]{MaMa}
%%\S \ref{sse-other}
where computing the drift is elementary and does not require any information
on the harmonic measure. 
According to \cite[Eq. (21)]{MaMa},
%%\eref{eq-elementary-general}
we have:
\begin{equation}\label{eq-elementary-z2z3}
\gamma = \frac{2p(1-2p)}{2-p} \:. 
\end{equation}

%%Here, computing the drift is elementary. Denote by $A$, resp. $B$,  the set of elements
%%of $\Z/2\Z\star \Z/3\Z$ whose normal form ends with $a$, resp. with
%%$b$ or $b^2$. When we are far from the unit element 1, the random walk on
%%$\Z/2\Z\star \Z/3\Z$ induces a Markov chain on $\{A,B\}$ with
%%transition matrix
%%\[
%%P=\left[ \begin{array}{cc} 0 & 1 \\ 1-p & p \end{array}\right]\:.
%%\]
%%For instance we have $P_{B,A}=1-p$, meaning that $P\{X_{n+1}\in A \mid X_n \in
%%B\}=1-p$, where $(X_n)_n$ is a realization of the random walk on
%% $\Z/2\Z\star \Z/3\Z$. 
%%The stationary distribution of $P$ is $\pi=[(1-p)/(2-p),1/(2-p)]$, and by
%%the Ergodic Theorem for Markov Chains, we have $\lim_n [\ P\{X_n\in
%%  A\} \ ,\ P\{X_n\in B\} \ ]=\pi$. The value of the drift follows readily:
%%\begin{equation}\label{eq-elementary-z2z3}
%%\gamma = \lim_n \frac{1}{n} \sum_{i=0}^{n-1}
%%E[X_{i+1}-X_i] + \frac{X_0}{n} = \lim_n E[X_n-X_{n-1}]= \frac{2p(1-2p)}{2-p} \:. 
%%\end{equation}
%%In particular for the simple random walk defined by
%%$\mu(a)=\mu(b)=\mu(b^2)=1/3$, we obtain $\gamma=2/15$. This last value was
%%derived in \cite{NeVo} by a direct evaluation and an asymptotic
%%analysis of $P\{|X_n|=k\}$. 
%%As a curiosity, observe that the drift is maximized for $p=2-\sqrt{3}$
%%and is then equal to $14-8\sqrt{3}=0.143594\cdots$. 

To go further, let us compute explicitly the harmonic measure.
Consider the Traffic Equations \cite[Eq. (11)]{MaMa}. 
%%\eref{eq-traffic} 
By symmetry, we should look for
a solution $r$ satisfying $r(b)=r(b^2)$. With this simplication, the
Traffic Equations yield a single equation:
$r(a)=2(1-p)r(b)$. So the solution is:
\[
\bigl[r(a),r(b),r(b^2)\bigr] = \Bigl[ \ \frac{1-p}{2-p}\ , \
  \frac{1}{2(2-p)}\ , \ \frac{1}{2(2-p)}\ 
  \Bigr]\:.
\]
Observe that the transition matrix defining the Markovian
multiplicative harmonic measure is (the elements of $\Sigma$ being
ranked in the order $\{a,b,b^2\}$):
\[
P=\left[ \begin{array}{ccc} 0 & 1/2 & 1/2 \\ 1 & 0 & 0 \\ 1 & 0 & 0
  \end{array} \right] \:.
\]
Hence the harmonic measure is the uniform measure except for the 
distribution of the initial element of the normal form. Observe that
the harmonic measure is never stationary. Indeed, the stationary
distribution of $P$ is $\pi=[1/2,1/4,1/4]$, and the equation $r=\pi$
has no solution in $p$. 

\medskip

Now let us consider the general nearest neighbor random walk on
$\Z/2\Z\star \Z/3\Z$ defined by $\mu(a)=1-p-q, \mu(b)=p,\mu(b^2)=q$. 
The direct computation of the drift is not possible
anymore. However, \cite[Theorem 3.5]{MaMa}
%%Theorem \ref{th-jackpot} 
applies and, furthermore, 
the harmonic measure can be explicitly
determined. 
Then we derive the drift using \cite[Corollary 3.6]{MaMa}. 
%%Corollary \ref{co-drift} 

\begin{figure}[ht]
\[ \epsfxsize=190pt \epsfbox{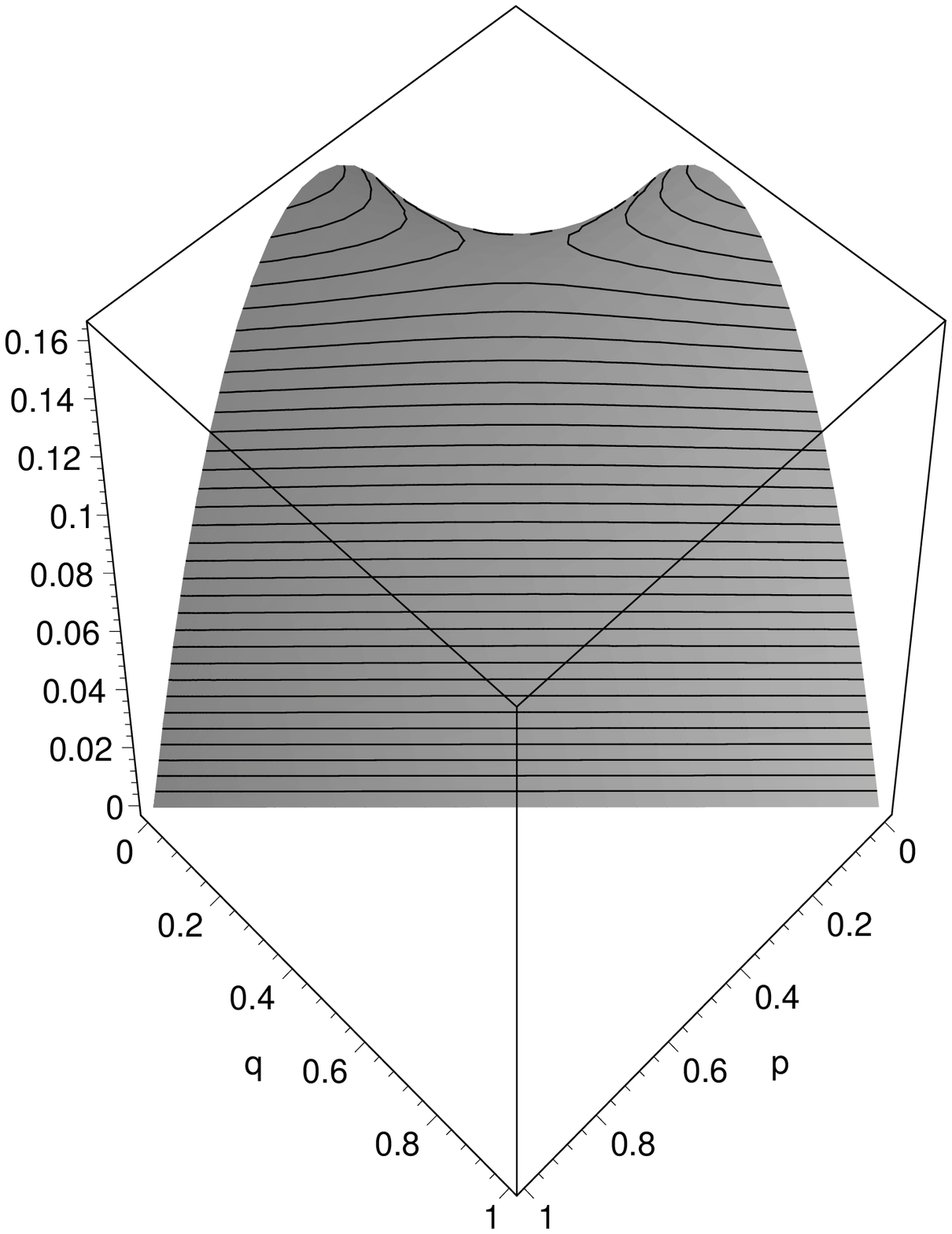} \qquad
 \epsfxsize=190pt \epsfbox{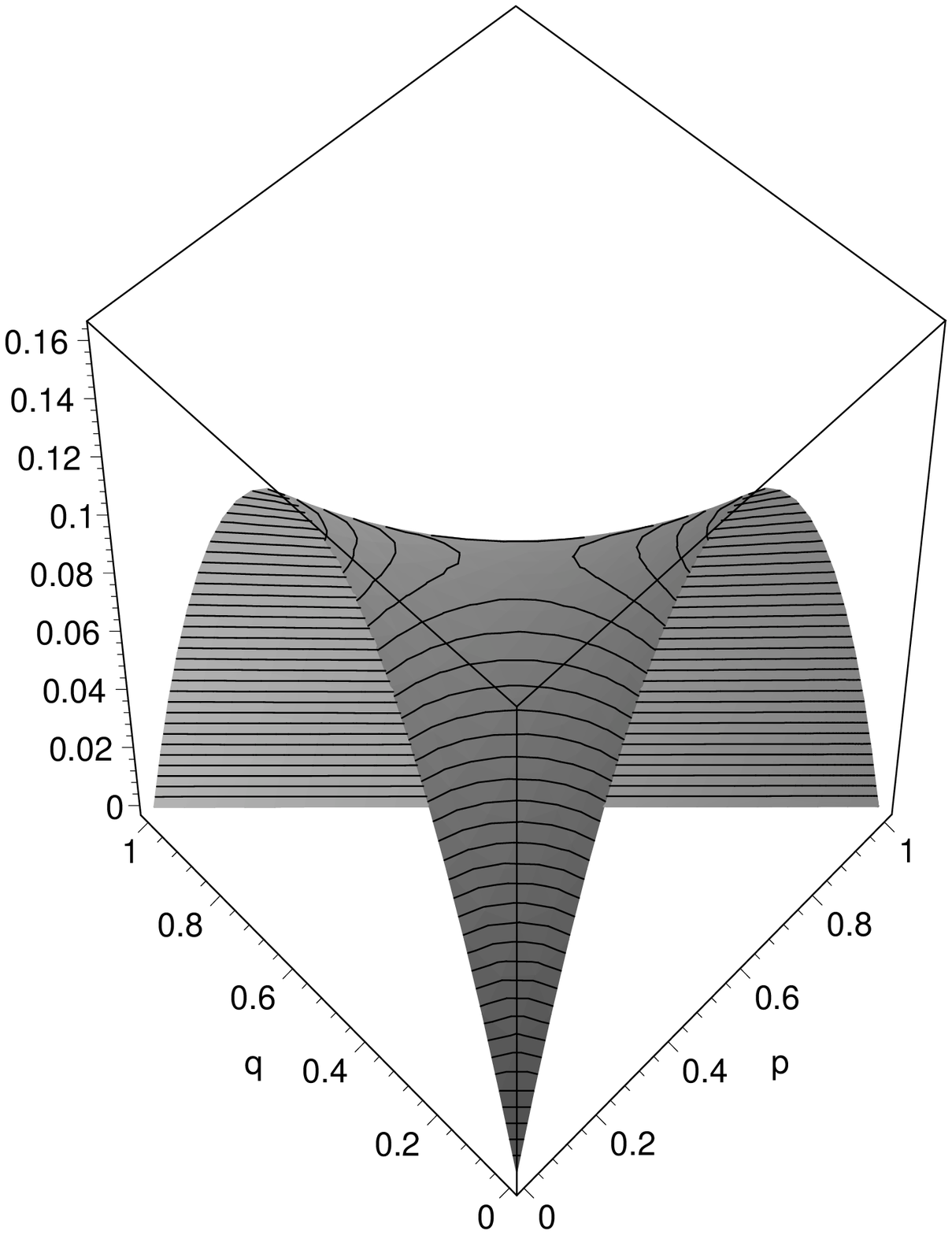} \]
\caption{The drift of the random walk $(\Z/2\Z\star \Z/3\Z,\mu)$ as a
  function of $p=\mu(b)$ and $q=\mu(b^2)$, from two different angles.}
\label{fi-z2z3drift} 
\end{figure}

%% It will enable us to get a closed form expression for the drift using
%% Proposition \ref{pr-drift}. 
The Traffic Equations \cite[Eq. (11)]{MaMa}
%%\eref{eq-traffic} 
are:
\begin{eqnarray}\label{eq-consz2z3}
r(a) & = & (1-p-q) (r(b)+r(b^2)) + q r(b) + p r(b^2) \nonumber \\
r(b) & = & pr(a) + q r(b^2) + (1-p-q) \frac{r(a)r(b)}{r(b)+r(b^2)} \\
r(b^2) & = & qr(a) + p r(b) + (1-p-q) \frac{r(a)r(b^2)}{r(b)+r(b^2)} \:.\nonumber
\end{eqnarray}
Substituting $1-r(b)-r(b^2)$ for $r(a)$ in the first equation, we get
$r(b^2)=(1-(2-p)r(b))/(2-q)$. Plugging this into the second (or third)
equation  allows us to compute $r(b)$. After some simplifications, we get:
\begin{eqnarray*}
r(a) & = & \frac{p^2+q^2-2pq-p-q+4 - \Delta_1}{2\Delta_2} \\
r(b) & = & \frac{q^3-3q^2 +p^2q -5pq+2p+6q -
  (2-q)\Delta_1}{2(q-p)\Delta_2} \\
r(b^2) & = & \frac{p^3 -3p^2+pq^2 -5pq +6p + 2q -
  (2-p)\Delta_1}{2(p-q)\Delta_2}  \:,
\end{eqnarray*}
with
\begin{eqnarray*}
\Delta_1 & = &
\sqrt{p^4+q^4-2p^3-2q^3+2p^2q^2-6p^2q-6pq^2+5p^2+5q^2+6pq} \\
\Delta_2 & = & p^2+q^2-pq-2p-2q+4 \:.
\end{eqnarray*}
Set $r=\mu(a)=1-p-q$. The drift is then:
\begin{equation}\label{eq-driftz2z3}
\gamma =  \frac{2r\Bigl( pq-p-q+
  \sqrt{(p^2+q^2)(3+(r+p)^2+(r+q)^2)+2pq(2r+1)}  
\Bigr)}{(r+p)^2+(r+q)^2 -pq+2  } \:.
\end{equation}
For instance, the drift is maximized for $r=z_0,p=1-z_0,q=0$ (or
$r=z_0,p=0,q=1-z_0)$), where $z_0$ is the root of
$[z^6+12z^4-4z^3+47z^2-48z+12]$ whose numerical value is
$0.490275\cdots$. The corresponding numerical value of the drift is
$\gamma_{\max}= 0.163379\cdots$. This was not a priori obvious!
%This may or may not be intuitive!

\par

We are in the domain of application of \cite[Proposition 3.9]{MaMa}.
%%Proposition \ref{pr-two} 
On the other hand, the harmonic measure is never shift-invariant.
This can be proved by showing directly that the Stationary
Traffic Equations have no solution.

\subsection{Random walks on $\Z/3\Z\star \Z/3\Z$}\label{sse-z3z3}

\begin{figure}[ht]
\[ \epsfxsize=240pt \epsfbox{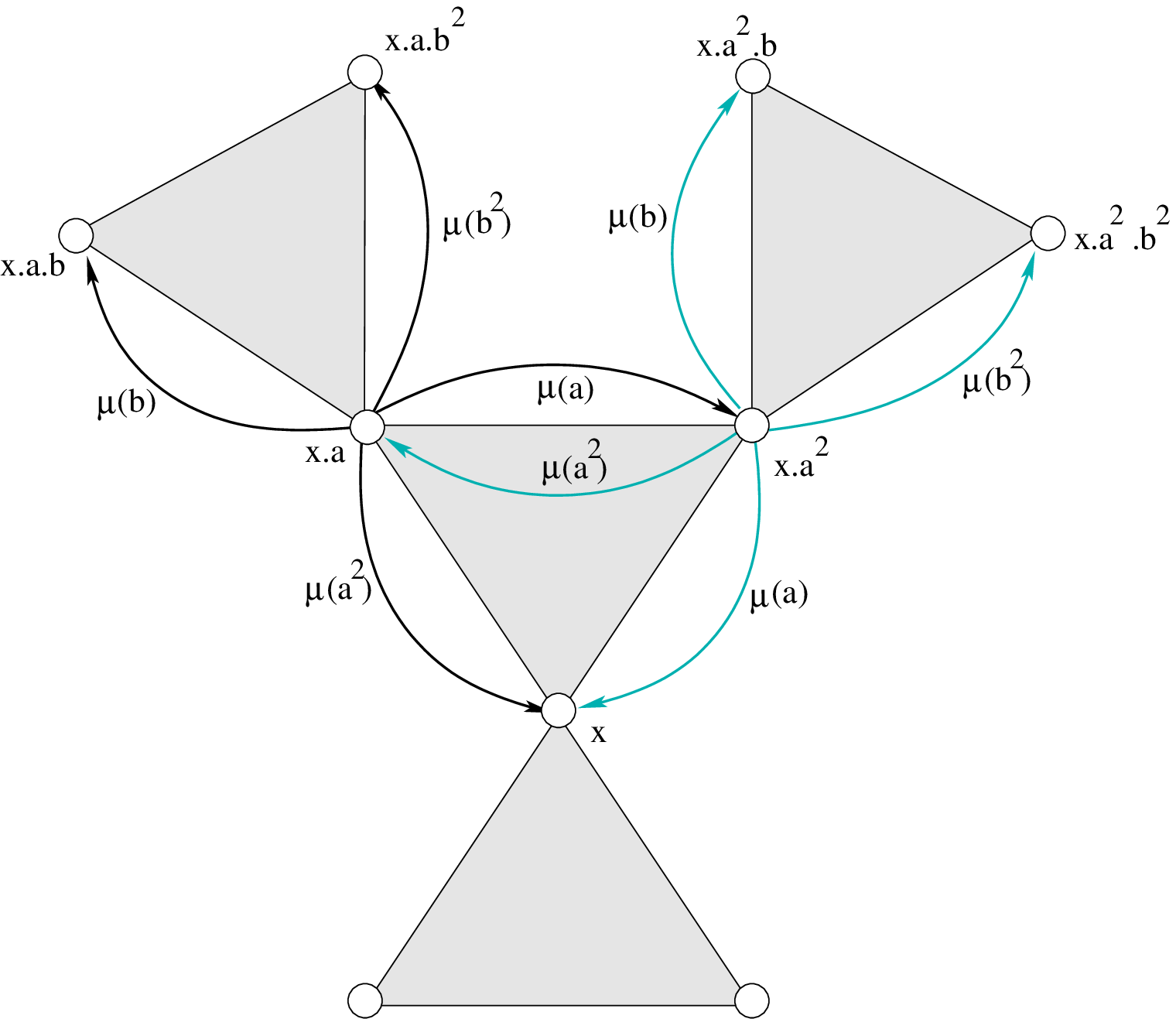} \]
\caption{A nearest neighbor random walk on $\Z/3\Z\star \Z/3\Z$.}
\label{fi-z3} 
\end{figure}

%% We turn our attention to a more complex case: the nearest neighbour
%% random walk $(\Z/3\Z\star \Z/3\Z,\mu)$ with the symmetry
%% $\mu(a)=\mu(b)=p,\mu(a^2)=\mu(b^2)=q=1/2-p$. Assume that $p\neq
%% q$. As opposed to the previous case, we cannot compute the drift
%% directly. Indeed we have
%% \begin{eqnarray*}
%% &&E[X_{i+1}-X_i\mid X_i=a] = E[X_{i+1}-X_i\mid X_i=b]=1/2-q \\
%% &&E[X_{i+1}-X_i\mid X_i=a^2] = E[X_{i+1}-X_i\mid X_i=b^2]=1/2-p\:,
%% \end{eqnarray*}
%% see Figure \ref{fi-z3}. However we are still able to determine the
%% harmonic measure explicitly. Then we will derive the drift using
%% Proposition \ref{pr-drift}.

Consider the probability $\mu$ such that
$\mu(a)=\mu(b)=p,\mu(a^2)=\mu(b^2)=q=1/2-p$. According to 
%%\ref{pr-stat}
\cite[Proposition 3.8]{MaMa}, the harmonic measure is ergodic.
Let us solve the Stationary Traffic Equations \cite[Eq. (12)]{MaMa}.
%%\eref{eq-trafficstat}
By obvious symmetry, we should have $r(a)=r(b)$ and $r(a^2)=r(b^2)$. Set
$r(i)=r(a^i)+r(b^i)$ for $i=1,2$. Using that $r(2)=1-r(1)$, the
Stationary Traffic Equations yield a single equation:
%% The  Stationary Traffic Equations become
%% \[
%% \begin{cases}
%% r(1) & = \ qr(2) + p + r(1) (qr(1)+pr(2)) \\
%% r(2) & = \ pr(1) + q + r(2) (qr(1)+pr(2))
%% \end{cases}\:.
%% \]
%% We rewrite the equations using that $r(2)=1-r(1)$ and $q=1/2-p$. We
%% get a single equation:
\[
\bigl(2p-\frac{1}{2}\bigr) r(1)^2 + \bigl(  2p -\frac{3}{2} \bigr)r(1)
-\frac{1}{2} =0\:.
\]
The unique solution lying between 0 and 1 gives
\begin{equation}\label{eq-z3z3a=b}
r(a)=r(b)=\frac{4p-3+\sqrt{16p^2-8p+5}}{4(4p-1)}, \quad
r(a^2)=r(b^2)=\frac{4p+1-\sqrt{16p^2-8p+5}}{4(4p-1)}\:.
\end{equation}
The harmonic measure is the ergodic Markovian multiplicative measure
associated with $r$. In particular, the drift is
\begin{equation}\label{eq-driftz3z3a=b}
\gamma \ = \ \frac{1}{2}-p + \bigl(2p-\frac{1}{2}\bigr)r(1) \ = \ -\frac{1}{4} +
\frac{1}{4}\sqrt{16p^2-8p+5}\:.
\end{equation}

\medskip

Now consider the nearest neighbour
random walk $(\Z/3\Z\star \Z/3\Z,\mu)$ with another symmetry:
$\mu(a)=\mu(a^2)=p,\mu(b)=\mu(b^2)=q=1/2-p$.
Here, it is not difficult
to check that there is no solution to the Stationary Traffic Equations \cite[Eq. (12)]{MaMa}.
%%\eref{eq-trafficstat}.
Let us turn our attention to the general Traffic Equations \cite[Eq. (11)]{MaMa}.
%%\eref{eq-traffic}. 
By obvious symmetry, we should look for a solution $r$ satisfying $r(a)=r(a^2)$ and
$r(b)=r(b^2)$. It implies that $r(\Sigma\moins \Sigma_a) =
r(b)+r(b^2)=2r(b)$ and $r(\Sigma\moins \Sigma_b) =r(a)+r(a^2)=2r(a)$.
The Traffic Equations get simplified as follows:
\[
r(a)  =   pr(a) + (2p+q) r(b) , \qquad
r(b)   =   qr(b) + (2q+p) r(a) \:.
\]
The solution to the above is:
\[
r(a)=r(a^2)=\frac{1+2p}{6}, \quad r(b)=r(b^2)=\frac{1-p}{3}\:.
\]
The harmonic measure is the Markovian multiplicative measure
associated with $r$.
We have $\mu^{\infty}(u_1u_2\cdots u_k \Sigma^{\N})= r(u_1)(1/2)^{k-1}$,
i.e. $\mu^{\infty}$ is the uniform measure except for the 
distribution of the initial element of the normal form. 
Using \cite[Eq. (19)]{MaMa},
%%\eref{eq-drift3}
the drift is
\begin{equation}\label{eq-elementary-z3z3}
\gamma = 4pq = 2p(1-2p)\:.
\end{equation}
The above formula for $\gamma$ could have been obtained using \cite[Eq. (21)]{MaMa}
%%\eref{eq-elementary-general} 
since $\mu$ is uniform on each of the two cyclic groups. 

\medskip

\begin{figure}[ht]
\[ \epsfxsize=240pt \epsfbox{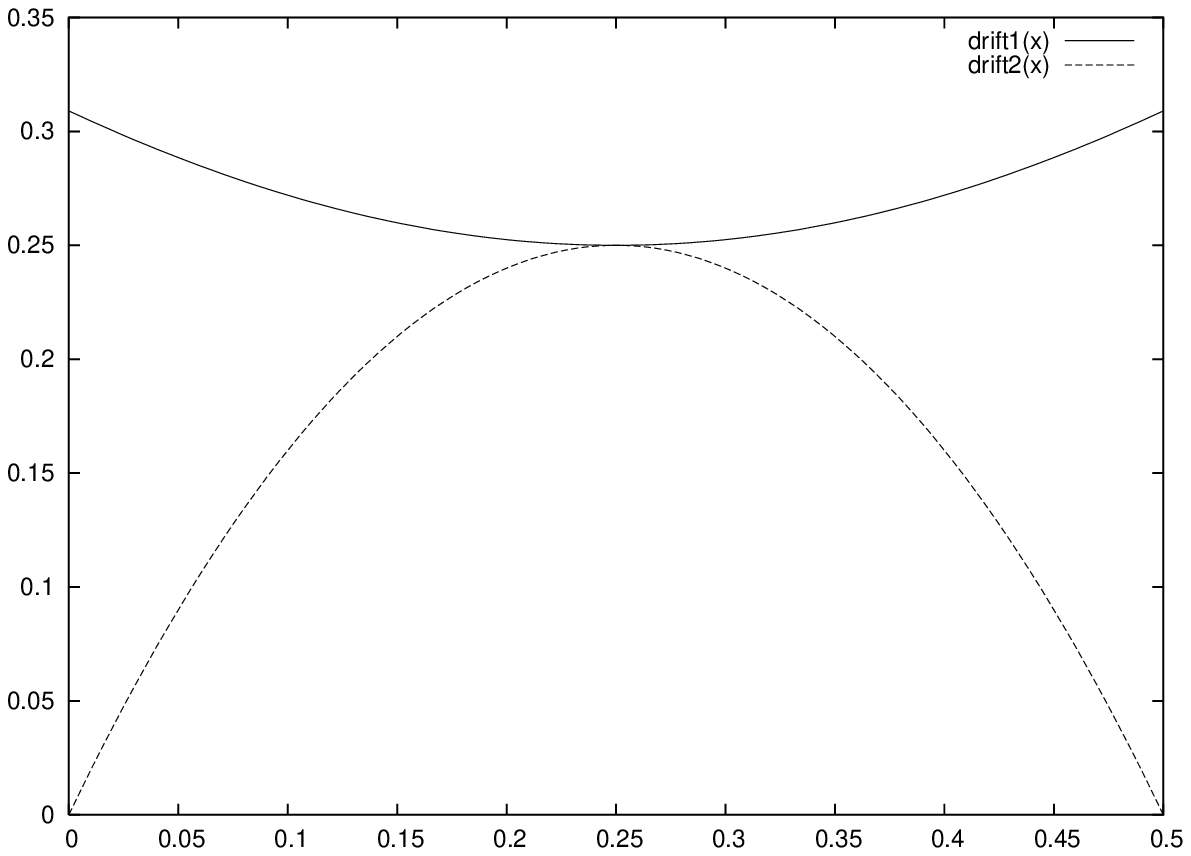} \]
\caption{The drift of $(\Z/3\Z\star \Z/3\Z,\mu)$ as a function of
  $p$.}\label{fi-driftz3} 
\end{figure}

We have represented in Figure \ref{fi-driftz3} the drift as a function
of $p$ for the two above cases: (i) $\mu(a)=\mu(b)=p,
\mu(a^2)=\mu(b^2)=1/2-p$ (drift1); (ii) $\mu(a)=\mu(a^2)=p,
\mu(b)=\mu(b^2)=1/2-p$ (drift2). 

At last, consider the case $\mu(a)=p,\mu(a^2)=q,$ and
$\mu(b)=\mu(b^2)=(1-p-q)/2$. Solving explicitly the Traffic Equations
is feasible but provides formula which are too lengthy to be
reproduced here. However, for the drift, several simplifications occur,
and we obtain the following formula: 
\[
\gamma= 2(1-p-q)\sqrt{\frac{p^2+q^2+pq}{p^2+q^2-2pq+3}}\:.
\]
The graph of $\gamma$ as a function of $p$ and $q$ looks like the one
in Figure \ref{fi-z2z3drift}.  

For the general nearest neighbor random walk on $\Z/3\Z\star \Z/3\Z$, we
did not succeed in solving completely the Traffic Equations. 

\subsection{The simple random walks on $\Z/k\Z\star \Z/k\Z$}
\label{sse-simplezn} 

\begin{figure}[ht]
\[ \epsfxsize=250pt \epsfbox{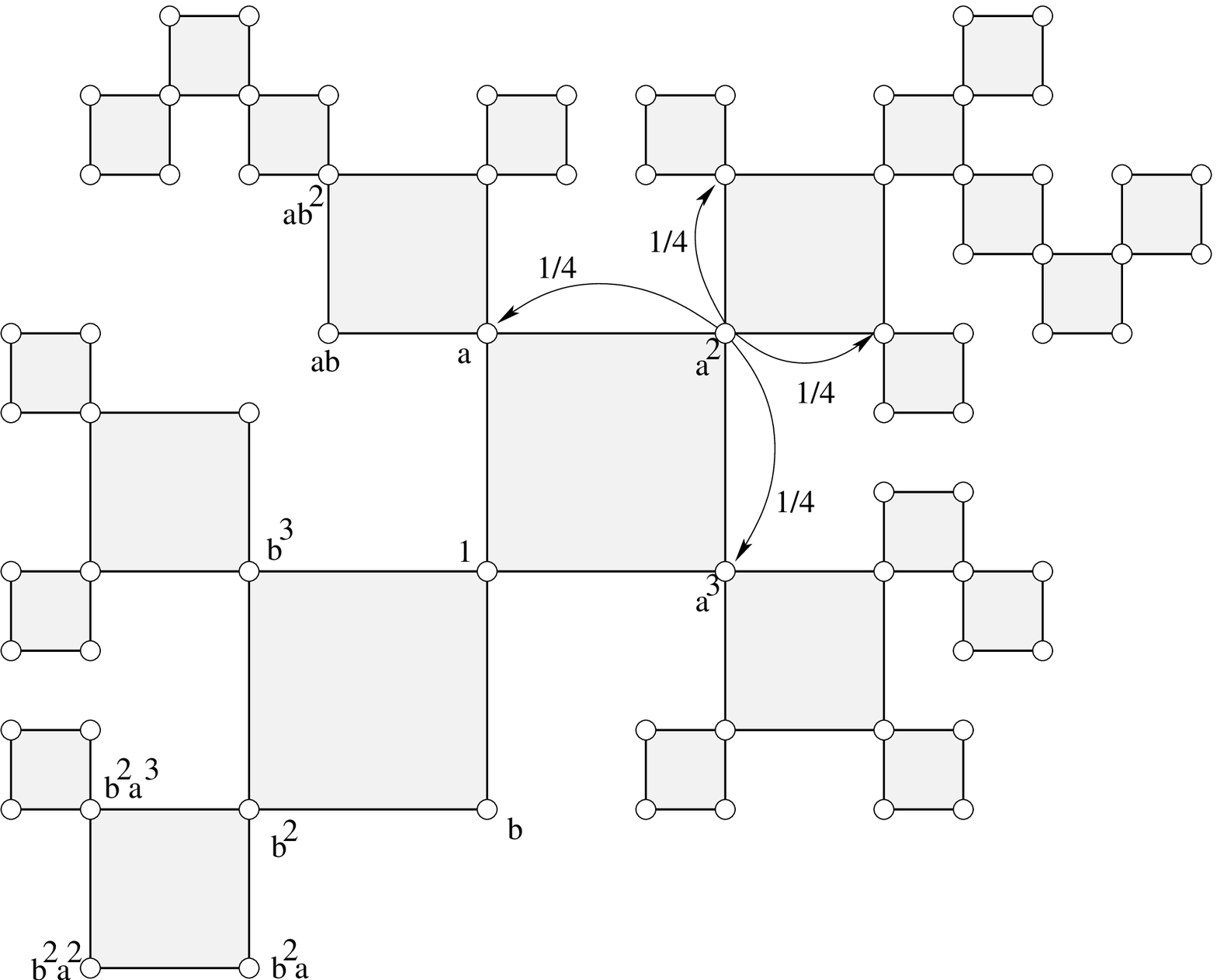} \]
\caption{The simple random walk on $\Z/4\Z\star \Z/4\Z$.}\label{fi-cayleyz4z4}
\end{figure}

We now consider the whole family of groups $\Z/k\Z\star \Z/k\Z, k\geq
3$. However, we have to compromise by
considering only simple random walks (with respect to a minimal set of
generators). In  Figure \ref{fi-cayleyz4z4}, we show this simple 
random walk in the case $\Z/4\Z\star \Z/4\Z$.
\par
Here we obtain ``semi-explicit'' formulas: we
define recursively a family of polynoms of one variable $(F_k)_k$, and 
the harmonic measure is expressed as a function of the 
unique solution in $(0,1)$ of $F_k=1$. For $k\geq 6$, we have no
closed form formula for this unique root. 

%% These results are central
%% in the study of Artin groups over 2 generators, see \S \ref{se-artin}. 

\medskip

Consider the free product $G_1\star G_2=\Z/k\Z\star \Z/k\Z$. Set
$\Sigma_1=G_1\moins\{1_{G_1}\}=\{a,\ldots ,a^{k-1}\}$,
  $\Sigma_2=G_2\moins\{1_{G_2}\}=\{b,\ldots ,b^{k-1}\}$, and
    $\Sigma=\Sigma_1\sqcup \Sigma_2$. Consider the simple random walk
    $(\Z/k\Z\star \Z/k\Z,\mu)$ with
    $\mu(a)=\mu(b)=\mu(a^{-1})=\mu(b^{-1})=1/4$. 

%% See Figure \ref{fi-cayleyz4z4}. 

%% We consider the simple random walk
%%     $(\Z/k\Z\star \Z/k\Z,\mu)$ with
%%     $\mu(a)=\mu(b)=\mu(a^{-1})=\mu(b^{-1})=1/4$. 
%% In Figure \ref{fi-cayleyz4z4}, we have represented the Cayley graph of
%% $\Z/4\Z\star
%% \Z/4\Z$ with respect to the generators $\{a,a^{-1},b,b^{-1}\}$. At each step, the random walk
%% moves to each one of the 4 neighboring vertices with probability 1/4.

According to \cite[Proposition 3.8]{MaMa},
%%Proposition \ref{pr-stat} 
the harmonic measure is ergodic. To determine it, we need to solve the simpler Stationary
Traffic Equations \cite[Eq. (12)]{MaMa}.
%%\eref{eq-trafficstat}. 
Set $r(i)=r(a^i)+r(b^i)=2r(a^i)=2r(b^i)$. 
The Stationary Traffic Equations are: 
\begin{equation}\label{eq-master}
r= r \Bigl[ \frac{1}{4} 
\left[\begin{array}{ccccc}
0 & 1 &  & &   \\
1 & 0 & 1 &  &   \\
  & \ddots & \ddots & \ddots &   \\
&&1&0 & 1 \\
&&& 1 & 0 
\end{array}\right] + \frac{1}{4} 
\left[\begin{array}{ccccc}
1 & 0 & \cdots & 0 & 1 \\
1 & 0 & \cdots & 0 & 1 \\
\vdots & \vdots & & \vdots & \vdots \\
1 & 0 & \cdots & 0 & 1 \\
1 & 0 & \cdots & 0 & 1
\end{array}\right] + \frac{1}{4}
\left[\begin{array}{cccc}
r(1) & r(2) & \cdots & r(k-1) \\
0 & \cdots & \cdots  & 0 \\
\vdots  & &   & \vdots \\
0 & \cdots & \cdots  & 0 \\
r(1) & r(2) & \cdots & r(k-1) 
\end{array}\right] \Bigr]\:.
\end{equation}

If $[r(1),r(2),\ldots , r(k-1)]$ is a solution to \eref{eq-master}, then it is
clear that $[r(k-1),\ldots,r(2),r(1)]$ is also a solution. By
uniqueness, we deduce that  $r(i)=r(k-i)$ for all $i$. In
particular, we have $r(1)=r(k-1)$. 
For convenience, set $r(0)=r(k)=1$. Then the
Equations \eref{eq-master} rewrite as:
\begin{equation}\label{eq-aftsimp}
2(2-r(1))r  =  [ r(0)+
    r(2),r(1)+r(3),\ldots,r(k-3)+r(k-1) , r(k-2)+r(k)]\:,
\end{equation}
that is, for all $i\in
\{1,k-1\}$,
\begin{equation}\label{eq-sym}
r(i+1)+r(i-1)=2(2-r(1))r(i) \:.
\end{equation}

Consider the applications $F_n: [0,1]\rightarrow \R, n\in \N,$ defined by 
\begin{equation}\label{eq-rec}
F_0(x)=1, \quad F_1(x) =x, \quad \forall n\geq 2, \ F_n(x) = 2(2-x)F_{n-1}(x)
-F_{n-2}(x)\:.
\end{equation}

To illustrate, here are the first values of $F_i$:
\begin{eqnarray*}
F_2(x) & = & - 2 x^2 +4 x - 1 \\
F_3(x) & = & 4x^3 -16x^2 + 17 x -4 \\
F_4(x) & = & -8x^4 +48x^3 -96x^2 +72x-15 \\
F_5(x) & = & 16 x^5 -128x^4 + 380 x^3 -512x^2 +301x-56 \\
F_6(x) & = & -32x^6 + 320 x^5 - 1264 x^4 + 2496 x^3- 2554 x^2 + 1244 x - 209
\end{eqnarray*}

By construction, the unique solution $r\in \mathring{\cB}$ of the Traffic Equations
\eref{eq-sym} satisfies $r(i)=F_i(r(1))$ for all 
$i\in \{1,\ldots, k-1\}$. Hence, it is enough to determine $r(1)$. 
For convenience, set $x_k=r(1)$.

The equalities $r(i)=r(k-i)$ for all $i\in \{1,\ldots,k-1\}$ rewrite
as $F_i(x_k)=F_{k-i}(x_k)$ for all $i\in \{1,\ldots,k-1\}$. Now, using
\eref{eq-rec} twice, we get
\begin{eqnarray*}
F_k(x_k) & = & 2(2-x_k)F_{k-1}(x_k)
-F_{k-2}(x_k) \\
& = & 2(2-x_k)F_{1}(x_k)-F_{2}(x_k) \\
& = & 2(2-x_k)F_{1}(x_k)-\bigl[
  2(2-x_k)F_{1}(x_k) - F_0(x_k)\bigr] = F_0(x_k)=1\:.
\end{eqnarray*}

Next lemma shows that the equality $F_k(x_k)=1$ is actually a
characterization of $x_k$. 

\begin{lemma}\label{le-unique}
For $k\geq 3$, the equation $F_k(x)=1$ has a unique solution $x_k$ in $(0,1)$. 
\end{lemma}

\begin{proof}
Notice first that:
\begin{equation}
\label{eq-1/3}
F_k\bigl(\frac{1}{3}\bigr)=\frac{1}{3^k}, \ \ F_k(1)=1 \:.
\end{equation}
Now, let us prove that there exists a point
$x^k\in (1/3,1)$ such that $F_k$ is strictly increasing on $[0,x^k)$
  and strictly decreasing on $(x^k,1]$. In view of \eref{eq-1/3}, the
proof of the lemma will then follow directly. 

Consider the matrix 
\begin{equation}
A=\left[\begin{array}{cc} 2(2-x) & -1 \\ 1 & 0 \end{array}\right]\:,
\end{equation}
for $x\in [0,1]$. 
Rewriting \eref{eq-rec} in matrix form, we get,
\begin{equation}\label{eq-dir}
\left[\begin{array}{c} F_k \\ F_{k-1} \end{array}\right] = A
\left[\begin{array}{c} F_{k-1} \\ F_{k-2} \end{array}\right] = A^{k-1}
\left[\begin{array}{c} F_{1} \\ F_{0} \end{array}\right]\:.
\end{equation}
The eigenvalues of $A$ are 
\begin{equation}
\lambda_1  =  2-x + \sqrt{(1-x)(3-x)}, \quad \lambda_2 = 2-x - \sqrt{(1-x)(3-x)} \:.
\end{equation}
and $\lambda_1\neq \lambda_2$ for $x\neq 1$. By diagonalizing, we get
for $x\neq 1$,
\begin{equation}\label{eq-diag}
A= \frac{1}{\lambda_1-\lambda_2} \left[\begin{array}{cc} \lambda_1 &
  \lambda_2 \\ 1 & 1 \end{array}\right] \left[\begin{array}{cc} \lambda_1 &
  0 \\ 0 & \lambda_2 \end{array}\right] \left[\begin{array}{cc} 1 &
  -\lambda_2 \\ -1 & \lambda_1 \end{array}\right]\:.
\end{equation}
Using \eref{eq-dir} and \eref{eq-diag}, we obtain
\begin{eqnarray}\label{eq-fk}
F_k(x)&=& \left[\begin{array}{cc} \frac{1}{\lambda_1-\lambda_2} &
  \frac{1}{\lambda_1-\lambda_2} \end{array}\right]
  \left[\begin{array}{cc} \lambda_1^k & 
  0 \\ 0 & \lambda_2^k \end{array}\right] \left[\begin{array}{c}
  x-\lambda_2 \\
  \lambda_1-x\end{array}\right] 
=\frac{x-\lambda_2}{\lambda_1-\lambda_2} \lambda_1^k + 
 \frac{\lambda_1-x}{\lambda_1-\lambda_2}\lambda_2^k \:.
\end{eqnarray}

By differentiating \eref{eq-fk}, we get
\begin{equation}\label{eq-derived}
F_k'= \frac{\lambda_1^k
  \Bigl[(x-3+2(1-x)^{1/2}(3-x)^{1/2})k +2\Bigr] + \lambda_2^k \Bigl[
  (-x+3+2(1-x)^{1/2}(3-x)^{1/2})k +2\Bigr] }{2(3-x)^{3/2}(1-x)^{1/2}} \:.
\end{equation}
Set $f_1(x)=x-3+2(1-x)^{1/2}(3-x)^{1/2}$ and
$f_2(x)=-x+3+2(1-x)^{1/2}(3-x)^{1/2}$. 
It is easily checked that $f_1$ and $f_2$ are strictly decreasing on
$[0,1]$. Since $f_2(1)=2$, we have $f_2>0$ on $[0,1]$. On the other
hand $f_1(1/3)=0$, hence we have $f_1>0$ on $[0,1/3)$ and $f_1<0$ on
  $(1/3,1]$. Since we have $\lambda_1>0$ and $\lambda_2>0$ for all $x\in
[0,1]$, it follows immediately that $F_k'>0$ on
$[0,1/3]$ and has at most one zero in $[0,1]$. Now by differentiating
\eref{eq-rec}, we get that
\[
F'_k(x)=2(2-x)F'_{k-1}(x)-2F_{k-1}(x)-F'_{k-2}(x)\:.
\]
In particular $F'_k(1)=2F'_{k-1}(1)-2-F'_{k-2}(1)$. It follows easily
that for $k\geq 3$, we have $F'_k(1)<0$. 
Consequently, $F_k'$ has exactly one zero in $[0,1]$ which belongs to
$(1/3,1)$. 
\end{proof}

\begin{figure}[ht]
\[ \epsfxsize=220pt \epsfbox{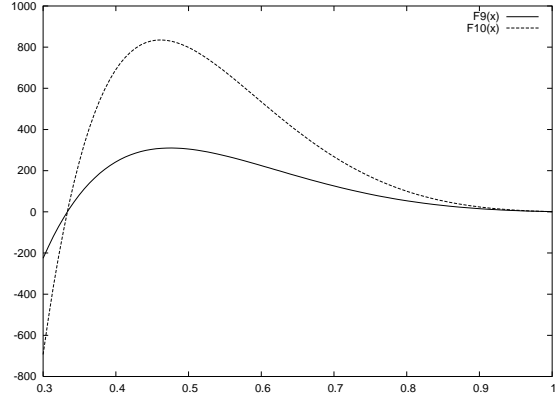} \]
\caption{The functions $F_9(x)$ and $F_{10}(x)$.}\label{fi-f9f10}
\end{figure}

To illustrate Lemma \ref{le-unique}, we have represented in Figure \ref{fi-f9f10}, the
functions $F_9$ and $F_{10}$. 

Let us collect the above results. The harmonic measure of
the simple random walk $(\Z/k\Z\star \Z/k\Z,\mu)$ is the ergodic
Markovian multiplicative measure associated with:
\[
r= [ x_k, F_2(x_k), F_3(x_k),\ldots , F_2(x_k), x_k ]
\]
where $x_k$ is the unique solution in $(0,1)$ of the equation
$F_k(x)=1$. 

Now let us compute the drift $\gamma_k$ using the harmonic measure,
see \cite[Corollary 3.8]{MaMa}.
%%Corollary \ref{co-drift}
We have
\begin{equation*}
\gamma_k = \sum_{a\in \Sigma} \frac{1}{4} \bigl[ -r(a^{-1})
  + \sum_{b\in \Sigma\moins
  \Sigma_a} r(b) \bigr] 
  =  4 \frac{1}{4} [ -\frac{x_k}{2} + \frac{1}{2} ] = 
  \frac{1-x_k}{2}\:.
\end{equation*}

\begin{lemma}
The drift $\gamma_k$ is a strictly increasing function of $k$ and satisfies
$\lim_k \gamma_k=1/3$. 
\end{lemma}

\begin{proof}
Fix $x\in (1/3,1)$. Recall that $F_n(x)>0$ for all $n$ (Lemma
\ref{le-unique}).    
Consider the expression
\eref{eq-fk}. We have $\lambda_1 >1$ and $0<\lambda_2<1$. It follows
easily that  $\lim_n F_n(x)=\lim_n \lambda_1^n
(x-\lambda_2)/(\lambda_1-\lambda_2)=+\infty$. 

Define 
$R_n=F_{n+1}(x)/F_{n}(x)$
for $n\in \N$. Using \eref{eq-rec}, we obtain 
\[
R_0= x, \quad \forall n\geq 1, \quad R_n=2(2-x) - \frac{1}{R_{n-1}}\:.
\]
Assume that there exists some $n$ such that $R_n\geq 1$. Then we have
$R_{n+1} \geq 2(2-x) -1 > 1$, since $x<1$. 
Since $\lim_n F_n(x)=+\infty$, we deduce that there exists $n^x \in \N\moins\{0\}$ such that
$F_n(x)$ is a strictly decreasing function of $n$ for $n< n^x$, and
a strictly increasing function of $n$ for $n> n^x$. 

Now, let us prove that $x_{k+1}<x_k$. We have 
$x_k\in (1/3,1)$, $F_1(x_k)=F_{k-1}(x_k)$. It implies that
$n^{x_k}\leq k-1$. In particular $F_{k+1}(x_k)>F_k(x_k)=1$. In view of
the behavior of $F_{k+1}$, as described in the proof of Lemma
\ref{le-unique}, it implies indeed that $x_{k+1}<x_k$. Therefore,
$\gamma_k$ is a strictly increasing function of $k$. 

Consider $x\in (1/3,1)$. We have seen that $\lim_n F_n(x)=\infty$. So
there exists $N_x$ such that $\forall n>N_x, F_n(x)>1$, and therefore
$\forall n>N_x, x_n < x$. It implies immediately that $\lim_n x_n =
1/3$ and $\lim_n \gamma_n =1/3$. 
\end{proof}

To summarize, we have proved the theorem below. 

\begin{theorem}\label{th-zkzk}
Consider the group $\Z/k\Z\star \Z/k\Z$, the generators of the two
cyclic groups being respectively $a$ and $b$. Consider the simple
random walk $(\Z/k\Z\star \Z/k\Z,\mu)$ with
$\mu(a)=\mu(a^{-1})=\mu(b)=\mu(b^{-1})=1/4$. 
The harmonic measure is
the Markovian multiplicative measure associated with $r$: $\forall i
\in \{1,\ldots, k-1\}, \ r(a^i)=r(b^i)=F_i(x_k)/2$, where $F_i$ is
defined in \eref{eq-rec} and $x_k$ is the unique solution in $(0,1)$
of $F_k(x)=1$. The drift is
$\gamma_k= (1-x_k)/2$. It is a strictly increasing function of $k$ and
$\lim_k \gamma_k=1/3$. 
\end{theorem}

To illustrate,
let us give the vector $r$ for $\Z/4\Z\star \Z/4\Z$:
\begin{equation}\label{eq-rz4z4}
\bigl[r(a),r(a^2),r(a^3)\bigr]=\bigl[r(b),r(b^2),r(b^3)\bigr]= \Bigl[
  \frac{3-\sqrt{5}}{8}, \frac{\sqrt{5}}{4} -\frac{1}{2}, \frac{3-\sqrt{5}}{8} \Bigr]\:.
\end{equation}
and for $\Z/5\Z\star \Z/5\Z$:
\[
\bigl[r(a),r(a^2),r(a^3),r(a^4)\bigr]=\bigl[r(b),r(b^2),r(b^3),r(b^4)\bigr]= \Bigl[
  \frac{5-\sqrt{13}}{8}, \frac{\sqrt{13}-3}{8}, \frac{\sqrt{13}-3}{8},
  \frac{5-\sqrt{13}}{8} \Bigr]\:.
\]
Now, here is a table of the first values of $\gamma_k$, given
either in closed form or numerically when no closed form could be
found. Set $\Z_k=\Z/k\Z$. 

\begin{center}
\begin{tabular}{|l|c|c|c|c|c|c|} \hline
 & $\Z_3\star \Z_3$ & $\Z_4\star \Z_4$ &   $\Z_5\star
 \Z_5$ &  $\Z_6\star \Z_6$ & $\Z_7\star \Z_7$ & $\Z_8\star
 \Z_8$ \\ \hline 
&&&&&& \\
$\gamma$ &  $1/4$ & $(\sqrt{5}-1)/4$  & $(\sqrt{13}-1)/8$  &
  0.330851... & 0.332515...  &  0.333062...\\ 
&&&&&& \\\hline
\end{tabular}
\end{center}

\paragraph{The limit case.} The limiting value of $1/3$ for $\gamma_k$
 as $k\rightarrow \infty$ was to be expected. Let us show why. The
`limit' of $\Z/k\Z \star \Z/k\Z$ as $k\rightarrow \infty$ is
 $\Z\star \Z$. The group $\Z\star \Z$ is the free group over two
 generators, its Cayley graph is $\T_4$, the homogeneous tree of
 degree 4. 
Analysing the simple random walk on $\Z\star \Z$ is
 straightforward. Let $a$ and $b$ be the respective generators of the two
 groups and let $\mu^{\infty}$ be the harmonic measure with respect to
 the set of generators $\Sigma=\{a,a^{-1},b,b^{-1}\}$. (The existence
 of $\mu^{\infty}$ is classical, see the discussion following \cite[Proposition 2.1]{MaMa}.)
%% Proposition \ref{pr-harm}
A word $u \in \Sigma^*$ is a normal form (or reduced) word if $u=u_1\cdots u_k$
with $\forall i, u_{i+1}\neq u_i^{-1}$. For any normal
form word $u\in \Sigma^k$, we have, by symmetry,  $\mu^{\infty}(u\Sigma^{\N})=1/4\times
1/3^{k-1}$. In particular the drift is $\gamma = 1/2$. 

However, our computations in $\Z/k\Z \star \Z/k\Z$
 correspond to the set of generators $\{a,\ldots, a^{k-1},
 b,\ldots,b^{k-1}\}$, and not to the set of generators
 $\{a,a^{-1},b,b^{-1}\}$. To get the analog in $\Z\star \Z$, we have
 to compute the drift with respect to the infinite set of generators
 $S=\{a^i,b^i,i\in \Z\moins\{0\}\}$. 
Let $X_{\infty}=x_1\cdots x_n\cdots$ be a r.v. distributed according
 to $\mu^{\infty}$. Let $\tau=\max\{ k \mid x_1=\cdots = x_k \}$. The
expectation $E[\tau]$ is the rate of reduction of the length when
 recoding a reduced word over $\Sigma$ in terms of $S$. 
Therefore the drift with respect to $S$ is equal to $(1/2\times 1/E[\tau])$. 
Using the above formula for $\mu^{\infty}$, we get $P\{\tau=k\}=
1/3^{k-1}\times 2/3$. It implies that $E[\tau]=3/2$. So 
the drift of the simple
random walk on $\Z\star \Z$ with respect to $S$ is $(1/3=1/2\times
2/3)$. 

\subsection{The simple random walks on $\Z/2\Z \star \Z/k\Z$}
\label{sse-z2zk}

We study simple (for a minimal set of generators) random walks on
$\Z/2\Z \star \Z/k\Z$. The results are ``semi-explicit'' as in \S \ref{sse-simplezn}.

The groups
$\Z/2\Z\star \Z/k\Z$ are known as the Hecke groups. A possible
faithful representation of $\Z/2\Z\star \Z/k\Z$, generalizing the one given in
\eref{eq-modular}, is
\[
a=\left[ \begin{array}{cc} 0  & -1 \\ 1 & 0 \end{array}\right], \qquad
b=\left[ \begin{array}{cc} 2\cos(\pi/k) & -1 \\ 1 & 0 \end{array}\right] \:,
\]
where $a$ and $b$ are the generators of $\Z/2\Z$ and $\Z/k\Z$
respectively.

\medskip

Consider the free product $\Z/2\Z\star \Z/k\Z$, with $a$ and $b$ the
respective generators of the two cyclic groups.  We consider the simple random walk
    $(\Z/2\Z\star \Z/k\Z,\mu)$ with
    $\mu(a)=\mu(b)=\mu(b^{-1})=1/3$. 

\begin{figure}[ht]
\[ \epsfxsize=200pt \epsfbox{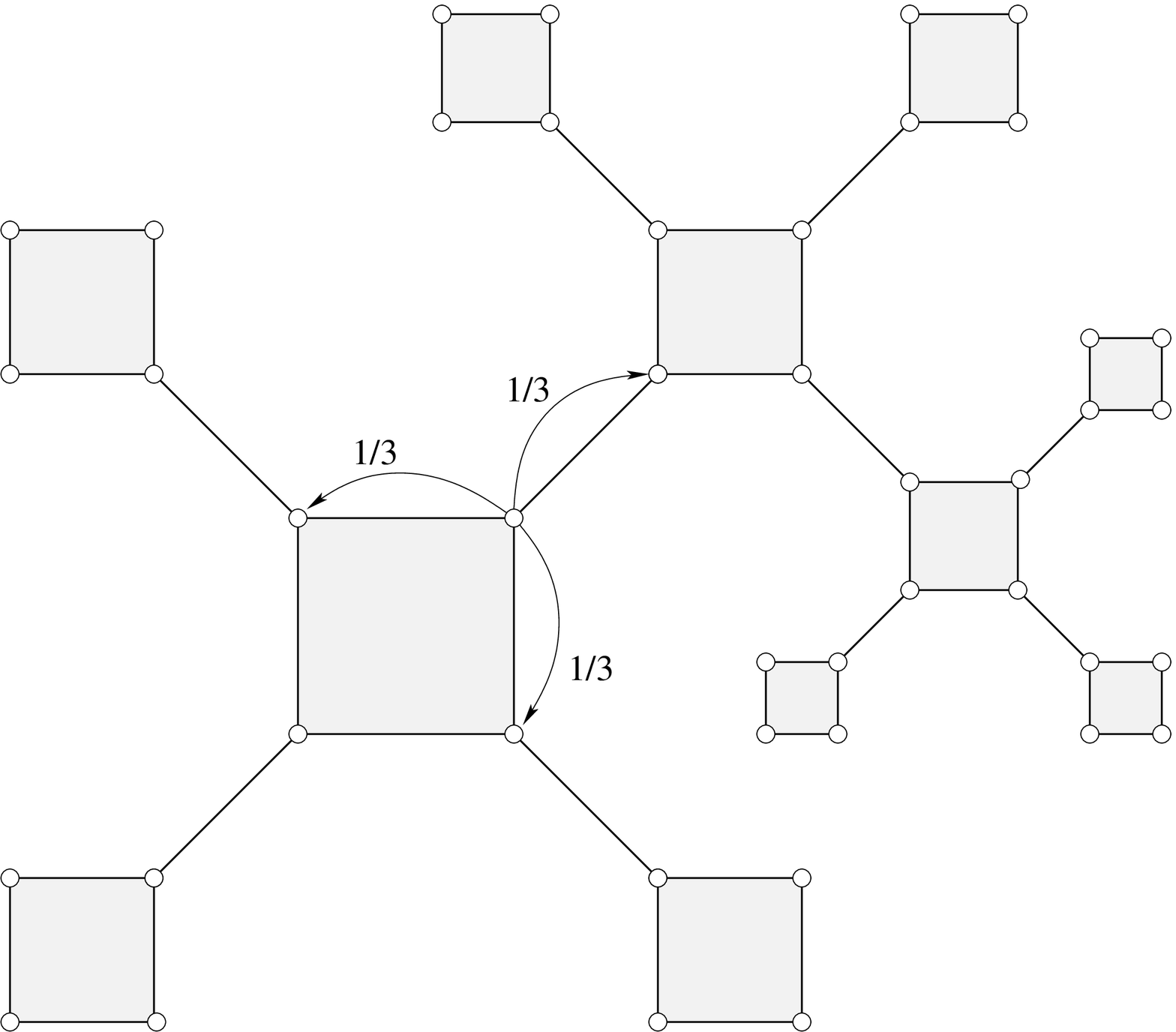} \]
\caption{The simple random walk on $\Z/2\Z\star \Z/4\Z$.}\label{fi-cayleyz2z4}
\end{figure}

We have seen above that even for $\Z/2\Z\star \Z/3\Z$, the harmonic
measure is not ergodic. Hence we should try to solve the general
Traffic Equations \cite[Eq. (11)]{MaMa}:
%%\eref{eq-traffic}
\begin{eqnarray*}
r(a) & = & \frac{1}{3}\Bigl[1-r(a) +  r(b) + 
r(b^{k-1})\Bigr]\\ 
r(b) & = & \frac{1}{3}\Bigl[ r(a)+r(b^2) +  r(b) \frac{r(a)}{1-r(a)}\Bigr]
\\
r(b^i) & = & \frac{1}{3}\Bigl[ r(b^{i-1})+r(b^{i+1}) +
r(b^i) \frac{r(a)}{1-r(a)}\Bigr], \qquad \forall i
=2,\ldots, k-2 \\ 
r(b^{k-1}) & = & \frac{1}{3}\Bigl[ r(b^{k-2})+r(a) +
r(b^{k-1}) \frac{r(a)}{1-r(a)}\Bigr] \:.
\end{eqnarray*}
%% Assume that $r=[r(a),r(b),r(b^2),\ldots ,r(b^{k-1}]$ is some solution to the
%% above with $r\in \cB$. Then, by symmetry of the equations, $[r(a),r(b^{k-1}),\ldots
%%   ,r(b^2),r(b)]$ is also a solution. But we know that there is a
%% unique solution to the Traffic Equations in $\cB$. We conclude
By symmetry, we must have $r(b^i)=r(b^{k-i})$ for all $i$. 
Hence the first equation gives: $r(a)= 1/4 + r(b)/2$. 
Now set $r(b^0)=r(b^k)=r(a)=1/4 + r(b)/2$. Rewriting the Traffic
Equations (except the first one), we get for $i
=1,\ldots, k-1$,
\begin{equation}\label{eq-symz2zk}
\Bigl[\frac{3-2r(b)}{8(1-r(b))}\Bigr] r(b^i)= r(b^{i-1})+r(b^{i+1}) \:. 
\end{equation}
This should be compared with \eref{eq-sym}. 
Consider the applications $G_n: [0,1]\rightarrow \R, n\in \N,$ defined by 
\begin{equation}\label{eq-recz2zk}
G_0(x)=\frac{1}{4}+\frac{x}{2}, \quad G_1(x) =x, \quad \forall n\geq
2, \ G_n(x) = \frac{8(1-x)}{3-2x}G_{n-1}(x) 
-G_{n-2}(x)\:.
\end{equation}
In illustration, here are the first few values of $G_i$:
\begin{eqnarray*}
(2x-3)G_2(x) & = & 7 x^2 - 7 x  + 3/4 \\
(2x-3)^2G_3(x) & = & 52x^3-100x^2+53x -6 \\
(2x-3)^3G_4(x) & = & 1552x^4-4416x^3+4296x^2-1600x+165
\end{eqnarray*}
From now on, the arguments become very close to the ones in \S
\ref{sse-simplezn}. Hence, we choose not to detail them. 
Observe just that the value of $x$ playing a special role (similar to
the one of 1/3 in \S \ref{sse-simplezn}) is now 1/6 with
$G_i(1/6)=1/(3\times 2^i)$.

\medskip

Next theorem summarizes the results. 

\begin{theorem}\label{th-z2zk}
Consider the group $\Z/2\Z\star \Z/k\Z$, the generators of the two
cyclic groups being respectively $a$ and $b$. Consider the simple
random walk $(\Z/2\Z\star \Z/k\Z,\mu)$ with
$\mu(a)=\mu(b)=\mu(b^{-1})=1/3$. 
Let the functions $G_i$ be defined as in \eref{eq-recz2zk} and let
$y_k$ be the
unique solution in $(0,1/2)$ of $G_{k-1}(y_k)=y_k$. 
The harmonic measure is
the Markovian multiplicative measure associated with $r$: $r(a) =
G_0(y_k), \ \forall i
\in \{1,\ldots, k-1\}, \ r(b^i)=G_i(y_k)$. The drift is $\gamma_k=
(1-2y_k)/3$. It is a strictly increasing function of $k$ and $\lim_k
\gamma_k =2/9$. 
\end{theorem}

Here is the vector $r$ for $\Z/2\Z\star \Z/4\Z$:
\[
\bigl[r(a),r(b),r(b^2),r(b^3)\bigr]= \Bigl[
  \frac{7-\sqrt{7}}{12}, \frac{2}{3} - \frac{\sqrt{7}}{6},
  \frac{-11+5\sqrt{7}}{12} , \frac{2}{3} - \frac{\sqrt{7}}{6}\Bigr]\:. 
\]
and for $\Z/2\Z\star \Z/5\Z$:
\[
\bigl[r(a),r(b),r(b^2),r(b^3),r(b^4)\bigr]= \Bigl[
  \frac{21-\sqrt{61}}{38}, \frac{23}{38} - \frac{\sqrt{61}}{19},
  \frac{-29+5\sqrt{61}}{76}, \frac{-29+5\sqrt{61}}{76},  \frac{23}{38}
  - \frac{\sqrt{61}}{19} \Bigr]\:.
\]
In next table, we give the drift of
the simple random walk on $\Z/2\Z \star \Z/k\Z$ for the small values
of $k$. The drift is given in closed form when possible, otherwise
numerically. Set $\Z_k=\Z/k\Z$. 

\begin{center}
\begin{tabular}{|l|c|c|c|c|c|c|} \hline
 & $\Z_2\star \Z_3$ & $\Z_2\star \Z_4$ &   $\Z_2\star
 \Z_5$ &  $\Z_2\star \Z_6$ & $\Z_2\star \Z_7$ & $\Z_2\star
 \Z_8$ \\ \hline 
&&&&&& \\
$\gamma$ &  $2/15$ & $(\sqrt{7}-1)/9$  & $(2\sqrt{61}-4)/57$  &
  0.213412... & 0.217921...  &  0.220101...\\ 
&&&&&& \\\hline
\end{tabular}
\end{center}

\paragraph{The limit case.} As in \S \ref{sse-simplezn}, the result $\lim_k
 \gamma_k =2/9$ can be recovered using only elementary arguments. 
The `limit' as $k\rightarrow \infty$ of $\Z/2\Z \star \Z/k\Z$ is
$\Z/2\Z \star \Z$. Let $a$ and $b$ be the respective generators of
$\Z/2\Z$ and $\Z$. The Cayley graph of $\Z/2\Z \star \Z$ with respect
to $\Sigma=\{a,b,b^{-1}\}$ is $\T_3$, the homogeneous tree of degree 3. 

From there, we proceed as in \S \ref{sse-simplezn}. 
The drift of the simple random walk on $\T_3$ is
obviously $1/3$. The harmonic measure is $\mu^{\infty} (u_1\cdots
u_k\Sigma^{\N})=1/3\times 1/2^{k-1}$ for any normal form word $u_1\cdots u_k\in
\Sigma^k$. Set $S=\{a,b^i,i\in \Z\moins \{0\}\}$. 
Let $X_{\infty}=x_1\cdots x_n\cdots$ be a r.v. distributed according
 to $\mu^{\infty}$. Let $\tau=\max\{ k \mid x_1=\cdots = x_k \}$. If
 $x_1=a$, then $\tau=1$. On the other hand, $P\{\tau=k\mid x_1=b \text{ or }
 b^{-1}\}= 1/2^{k}$. So we have $E[\tau\mid x_1=a]=1$ and
 $E[\tau\mid x_1\in\{b,b^{-1}\}]=2$. In a normal form word over
 $\Sigma$, we have a strict alternance of a single $a$ with blocks of
 $b$ or $b^{-1}$. Hence the expected reduction in length when recoding
 a normal form word over $S$ is $3/2$. 
We conclude that the drift of  the simple
random walk on $\Z/2\Z\star \Z$ with respect to $S$ is $(2/9=1/3\times
2/3)$.

\section{Entropy and Extremal Generators}\label{se-vershik}

Here, we give the proofs of the results announced in \cite[\S
  5.1]{MaMa}, and recalled below.

\paragraph{Statement of the results.} $ $ \\

\hspace*{1cm} $\bullet$ A. Consider a free product of two finite
groups, $G_1\star G_2$. 
Consider the probability $\mu_p$ on $\Sigma$ such
that $\mu_p = p U_{\Sigma_1} + (1-p) U_{\Sigma_2}$,
where $U_{\Sigma_i}$ is the uniform distribution on $\Sigma_i$, and
where $p\in (0,1)$. Then we have
$h/(\gamma v)=1$, for the whole family of probabilities $\mu_p$. This
does not extend to the free product of more than 2 groups, see B,
nor to minimal sets of generators $S \subset \Sigma$, see C.

\par

\hspace*{1cm} $\bullet$ B. Consider the group $\Z/2\Z \star \Z/2\Z \star \Z/2\Z$,
and let $a,b,$ and $c$ be the non-identity elements of the three
cyclic groups. Consider the family of probability measures $\mu_p,
p\in (0,1/2),$ defined by: $\mu_p(a)=\mu_p(b)=p, \mu_p(c)=1-2p$. Among
them, the
only probability such that $h/(\gamma v)=1$ is $\mu_{1/3}$. 

\par

\hspace*{1cm} $\bullet$ C. Consider the group $\Z/2\Z \star \Z/4\Z$. Let $a$ and
$b$ be respective generators of the two cyclic groups. The
minimal set of generators
$S=\{a,b,b^{-1}\}$ is extremal but $h/(\gamma v)<1$ for
all $\mu$ in $\cS$.
Consider now the group $\Z/3\Z \star \Z/4\Z$. The minimal set
of generators $S=\{a,a^{-1},b,b^{-1}\}$ is extremal. Actually, the only
symmetric probability measure on $S$ for which $h=\gamma v$ is
$\mu \,=\, p\,(\delta_{a}+\delta_{a^-1})+(1/2-p)(\delta_{b}+\delta_{b^-1})$ 
with $p\,=\,0.432692\cdots$ is the middle root of the
polynomial $5x^3-13x^2+7x-1$. Observe that $\mu \in \cS$. 
Last, consider the group $\Z/4\Z \star \Z/4\Z$. The minimal set
of generators $S=\{a,a^{-1},b,b^{-1}\}$ is not extremal. Indeed, 
\[
Q(S)=  \frac{5+\sqrt{5}}{4}\frac{\log \bigl(1/2 + \sqrt{5}/2\bigr)}
{\log(1+\sqrt{2})} = 0.987686\cdots \:.
\]

\paragraph{Proof of A.} 

The volume is
$v=\log(|\Sigma_1||\Sigma_2|)/2$. Set $r(\Sigma_i)=\sum_{x\in \Sigma_i}
r(x)$ and $q(\Sigma_i)= \sum_{x\in \Sigma_i}
q(x)$. 
By symmetry, for all $x\in \Sigma_i$, we have
$r(x)=r(\Sigma_i)/|\Sigma_i|$ and $q(x)=q(\Sigma_i)/|\Sigma_i|$. 
Using \cite[Eq. (7)]{MaMa} and \cite[Eq. (36)]{MaMa}
%%\eref{eq-drift2} and \eref{eq-entro2}
 we obtain
\begin{eqnarray*}
\gamma & = & \mu(\Sigma_1) \bigl( r(\Sigma_2)-\frac{r(\Sigma_1)}{|\Sigma_1|}
\bigr) + \mu(\Sigma_2) \bigl( r(\Sigma_1)-\frac{r(\Sigma_2)}{|\Sigma_2|}
\bigr) \label{eq-drif3}  \\
h & = & -\mu(\Sigma_1)\log\bigl(\frac{|\Sigma_1|}{q(\Sigma_1)}\bigr)
\bigl(r(\Sigma_2)-\frac{r(\Sigma_1)}{|\Sigma_1|}\bigr) -
\mu(\Sigma_2)\log\bigl(\frac{|\Sigma_2|}{q(\Sigma_2)}\bigr) 
\bigl(r(\Sigma_1)-\frac{r(\Sigma_2)}{|\Sigma_2|}\bigr) \:. \label{eq-entropy3}
\end{eqnarray*}
To go further, let us solve the Traffic Equations. 
We have:
\begin{equation*}
r(\Sigma_1) =  \frac{\mu(\Sigma_1) + \mu(\Sigma_2)/|\Sigma_2|}{1+
  \mu(\Sigma_1)/|\Sigma_1| + \mu(\Sigma_2)/|\Sigma_2|}\:, \quad 
r(\Sigma_2) =  \frac{\mu(\Sigma_2) + \mu(\Sigma_1)/|\Sigma_1|}{1+
  \mu(\Sigma_1)/|\Sigma_1| + \mu(\Sigma_2)/|\Sigma_2|}\:.
\end{equation*}
It follows immediately that:
\begin{eqnarray*}
&r(\Sigma_1)\bigl( \mu(\Sigma_2)+ \mu(\Sigma_1)/|\Sigma_1|\bigr) & = \
r(\Sigma_2)\bigl( \mu(\Sigma_1)+ \mu(\Sigma_2)/|\Sigma_2|\bigr) \\
\iff & \mu(\Sigma_1) \bigl( r(\Sigma_2) - r(\Sigma_1)/|\Sigma_1|
\bigr) & = \
\mu(\Sigma_2) \bigl( r(\Sigma_1) - r(\Sigma_2)/|\Sigma_2| \bigr) \:.
\end{eqnarray*}
Using the above equality, we can rewrite the drift and the entropy as
follows:
\begin{eqnarray*}
\gamma & = & 2\mu(\Sigma_1) \bigl( r(\Sigma_2)-r(\Sigma_1)/|\Sigma_1|
\bigr) \\
h & = & -\mu(\Sigma_1)(r(\Sigma_2)-r(\Sigma_1)/|\Sigma_1|)\bigl( \log(|\Sigma_1|/q(\Sigma_1))
+ \log(|\Sigma_2|/q(\Sigma_2)) \bigr) \\
& = & -\mu(\Sigma_1)(r(\Sigma_2)-r(\Sigma_1)/|\Sigma_1|) \log(|\Sigma_1||\Sigma_2|)\:,
\end{eqnarray*}
where we have used the identity $q(\Sigma_1)q(\Sigma_2)=1$ to obtain the
last equality. 
We conclude that:
$h/\gamma = \log(|\Sigma_1||\Sigma_2|)/2=v$. 

\medskip

\paragraph{Proof of B.} 

The volume of $\Z/2\Z\star \Z/2\Z \star \Z/2\Z$ is $v=\log(2)$. The
Traffic Equations can be explicitly solved and we obtain the drift
$\gamma_p$ and the entropy $h_p$ associated with the measure $\mu_p$:
\begin{eqnarray*}
\gamma_p & = & 3p^2 +\frac{3-9p}{2} + \frac{3p-1}{2}\sqrt{4p^2-20p+9} \\
h_p & = & \gamma_p\log(2) + \log\bigl( 3-2p+
  \sqrt{4p^2-20p+9}\bigr)\frac{p}{2}(-1+2p+\sqrt{4p^2-20p+9}) + \\
&& \qquad  \log\Bigl(-1 +
  \frac{\sqrt{4p^2-20p+9}}{1-2p}\Bigr) \frac{1-2p}{2}(3-2p-\sqrt{4p^2-20p+9})\:.
\end{eqnarray*}
Let us write $h_p=\gamma_p\log(2) + f(p)$. A direct analysis shows that
$f(p)=0$ if and only if $p=1/3$. We deduce that $h_p/\gamma_p=v$ if and
only if $p=1/3$. 
%% \begin{figure}[ht]
%% \[ \epsfxsize=200pt \epsfbox{z2z2z2rat.eps} \]
%% \caption{The ratio $h_p/(\gamma_p v)$  as a
%%   function of $p$ for $\Z/2\Z\star \Z/2\Z \star
%%   \Z/2\Z$.}\label{fi-z2z2z2ratio} 
%% \end{figure}

\medskip

\paragraph{Proof of C.} 

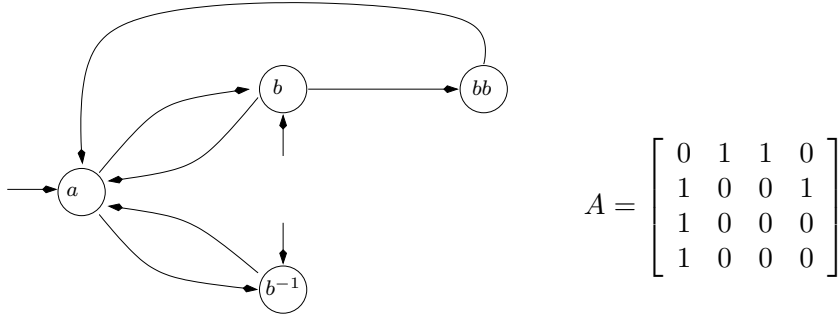
\begin{figure}[ht]
\[ \input{volume0.pstex_t} \]
\caption{An automaton enumerating the group elements according to
  $|\cdot|_S$. }\label{fi-volume} 
\end{figure}

Consider first the group $\Z/2\Z\star \Z/4\Z$ and let $a$ and $b$ be
the respective generators of the cyclic groups. Consider the minimal
set of generators $S=\{a,b,b^{-1}\}$. 
First of all, let us compute the volume $v(S)$ with respect to $S$. 
Deciding to write $b^2$ as $bb$, we enumerate the group elements using
the automaton of Figure \ref{fi-volume}. More precisely,
each path of length $n$ (starting with $a$, $b$, or $b^{-1}$)
corresponds to a different 
group element of length $n+1$ with respect to $S$.  
It follows that $v(S)=\log(\rho)$ where $\rho$ is the spectral radius of
$A$, the matrix associated with the automaton. We deduce that
$v(S)=\log(1/2+\sqrt{5}/2)$.

\begin{figure}[ht]
\[ \epsfxsize=200pt \epsfbox{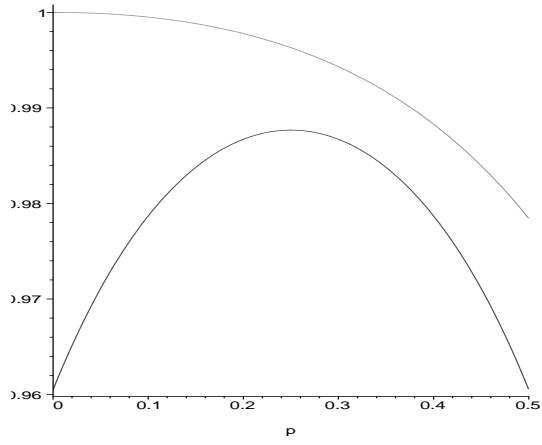} \]
\caption{The ratio $h(\mu)/(\gamma(S,\mu) v(S))$ for $(\Z/2\Z\star
  \Z/4\Z,\mu)$ (top) and $(\Z/4\Z\star
  \Z/4\Z,\mu)$ (bottom) as a function of $p=\mu(b)=\mu(b^{-1})$.}
\label{fi-zz} 
\end{figure}

Consider $\mu\in \cS$ with $\mu(a)=1-2p, \mu(b)=\mu(b^{-1})=p, p\in (0,1/2)$. 
The solution to the Traffic Equations is:
\begin{eqnarray*}
r(a) & = & \frac{ 2 -p - p^{2} - p\sqrt{5-6p+p^2}}{4(1-p)} \\
r(b)=r(b^{-1}) & = &
\frac {3-p - \sqrt{5-6p+p^2}}{4} 
  \\
r(b^2) & = & \frac {-4 + 5p - p^{2} +(2-p)\sqrt{5-6p+p^2}}{4(1-p)} \:.
\end{eqnarray*}
The value of the drift follows
(when
computing the drift, one must remember that the element $b^2$ 
has length 2 with respect to $S$):
\begin{eqnarray*}
\gamma(S,\mu) & = &  \frac{p(2p+5)}{2} -
\frac{p(2p+1)\sqrt{5-6p+p^2}}{2(1-p)} \:.
\end{eqnarray*}
Knowing $r$, 
the explicit formula for $h(\mu)$ is computable, but quite long. 
One checks directly on the resulting expression that 
$h(\mu)/(\gamma(S,\mu) v(S)) <1$ for all $p\in (0,1/2)$, and that 
$\lim_{p\searrow 0} h(\mu)/(\gamma(S,\mu) v(S)) =1$. See Figure \ref{fi-zz}. 
Hence $Q(S)=1$, the set of generators $S$ is extremal. 
%% A Taylor expansion
%% around 0 provides the following:
%% $h/(\gamma v)=1 - 4p^2/(81\log(3)) + O(p^3)$. 
%% Hence $Q(\Sigma)=1$, the set of generators $S$ is extremal. 

\medskip

Consider now the group $\Z/3\Z\star \Z/4\Z$ and the minimal
set of generators $S=\{a,a^{-1},b,b^{-1}\}$. With the same type of
arguments as above, we obtain that the volume with respect to $S$ is 
$v(S)=\log(\rho)$ where $\rho=2.21432\cdots$ is the real root of the polynomial $x^3-4x-2$.

Consider the measure $\mu\in \cS$ defined by 
$\mu(a)=\mu(a^{-1})=p, \mu(b)=\mu(b^{-1})=1/2-p, p\in (0,1/2)$.
The solution to the Traffic Equations is:
\begin{eqnarray*}
r(a)=r(a^{-1}) & = & \frac{ 5 + 4p - 4p^2 - (1-2p)\sqrt{9 + 4p - 4p^2} }{8(2-p)(1+p)} \\
r(b)=r(b^{-1}) & = &\frac { 5 - 2p - \sqrt{9 + 4p - 4p^2} }{4(2-p)} \\
r(b^2) & = & \frac { -7 - 6p + 4p^2 + 3\sqrt{9 + 4p - 4p^2} }{4(2-p)(1+p)} \:.
\end{eqnarray*}

Computing the drift cautiously as above gives:
\begin{eqnarray*}
\gamma(S,\mu) & = & \frac{(1-2p)(12 + 11p - 10p^2 - (4-5p)\sqrt{9 + 4p - 4p^2})}{4(2-p)(1+p)}  \:.
\end{eqnarray*}

With the explicit formula for $h(\mu)$, one checks directly that the
maximum of $h(\mu)/(\gamma(S,\mu) v(S))$ is $1$ and is obtained for
$p\,=\,0.432692\cdots$ which is the middle root of the polynomial $5x^3-13x^2+7x-1$.

\medskip

Last, consider the group $\Z/4\Z\star \Z/4\Z$ and the minimal
set of generators $S=\{a,a^{-1},b,b^{-1}\}$. Again, with the same type of
arguments as above, we obtain that the volume with respect to $S$ is 
$v(S)=\log(1+\sqrt{2})$. 

Consider the measure defined by $\mu(a)=\mu(a^{-1})=p,
\mu(b)=\mu(b^{-1})=1/2-p$ for some $p\in (0,1/2)$. 
The solution to the Traffic Equations is obtained formally. Consider
the two equations:
\begin{eqnarray*}
(-4+8p^2)X^4 +
(4+8p^2)X^3+(3+4p)X^2 - (2+2p)X- (1+2p) & = & 0 \\ 
(-2-8p+8p^2)X^4+(6-8p+8p^2)X^3+(5-4p)X^2+(-3+2p)X+2p-2 & = & 0
\end{eqnarray*}
Denote respectively by $T$ and $U$ the largest real root of the two
equations. (The exact formulas for $T$ and $U$ are too long to be
reproduced here.) We have: 
\[
r(a)=r(a^{-1})= \frac{1}{2(1+U)}, \ r(a^2)= \frac{1}{2U(1+U)}, \
r(b)=r(b^{-1})= \frac{1}{2(1+T)}, \ r(b^2)= \frac{1}{2T(1+T)}\:.
\]
The value of the drift and the entropy follows readily.  We get:
\begin{eqnarray*}
\gamma(S,\mu) & = & \frac{1-4p+2U(1-2p)}{2U(1+U)} + \frac{-1+4p+4pT}{2T(1+T)}
\\
h(\mu) & = &  \log\Bigl( \frac{(U+1)(2T+1)}{T(T+1)} \Bigr)
\frac{p(2U^2-U-1)}{U(U+1)}  + \log(U)\frac{p(U-1)}{U(U+1)}   \\ 
&  & \quad   + \log \Bigl( \frac{(T+1)(2U+1)}{U(U+1)} \Bigr)
\frac{(1/2-p)(2T^2-T-1)}{T(T+1)} +
\log(T)\frac{(1/2-p)(T-1)}{T(T+1)} \:.
\end{eqnarray*}
The ratio $h(\mu)/(\gamma(S,\mu) v(S))$ as a function of $p$ has been represented in
Figure \ref{fi-zz}. It is maximized for $p=1/4$. This was to be
expected since the two groups in the free product are isomorphic. 
For $p=1/4$, one has: $r(a)=r(b)=r(a^{-1})=r(b^{-1})  =
(3-\sqrt{5})/4,  \ r(a^2)=r(b^2) =  -1 + \sqrt{5}/2$, see
\eref{eq-rz4z4}. Hence,  
$\gamma(S,\mu) = (3-\sqrt{5})/2$ and
$h(\mu)=(5-\sqrt{5})\log \bigl(1/2 + \sqrt{5}/2\bigr)/4$. It follows that:
\[
Q(S)=  \frac{5+\sqrt{5}}{4}\frac{\log \bigl(1/2 + \sqrt{5}/2\bigr)}
{\log(1+\sqrt{2})} = 0.987686\cdots \:.
\]
In particular, the generators $S$ are not extremal.

%%          \bibliographystyle{plain}
%%           \bibliography{mairesse}
%%           \end{document}

\end{document}

%% file: volume0.pstex_t
\begin{picture}(0,0)%
\includegraphics{volume0.pstex}%
\end{picture}%
\setlength{\unitlength}{2763sp}%
\begingroup\makeatletter\ifx\SetFigFont\undefined%
\gdef\SetFigFont#1#2#3#4#5{%
  \reset@font\fontsize{#1}{#2pt}%
  \fontfamily{#3}\fontseries{#4}\fontshape{#5}%
  \selectfont}%
\fi\endgroup%
\begin{picture}(7718,2808)(2314,-4780)
\put(4711,-2801){\makebox(0,0)[lb]{\smash{\SetFigFont{8}{9.6}{\rmdefault}{\mddefault}{\updefault}$b$}}}
\put(6499,-2811){\makebox(0,0)[lb]{\smash{\SetFigFont{8}{9.6}{\rmdefault}{\mddefault}{\updefault}$bb$}}}
\put(4644,-4611){\makebox(0,0)[lb]{\smash{\SetFigFont{8}{9.6}{\rmdefault}{\mddefault}{\updefault}$b^{-1}$}}}
\put(2861,-3721){\makebox(0,0)[lb]{\smash{\SetFigFont{8}{9.6}{\rmdefault}{\mddefault}{\updefault}$a$}}}
\put(7501,-3886){\makebox(0,0)[lb]{\smash{\SetFigFont{11}{13.2}{\rmdefault}{\mddefault}{\updefault}$A = \left[ \begin{array}{cccc}0& 1 & 1 & 0\\ 1 & 0 & 0 & 1 \\ 1 & 0 & 0 &0 \\ 1 & 0 & 0 & 0 \end{array}\right] $}}}
\end{picture}